\documentclass[12pt,reqno,a4paper]{amsart}
\title[Counterexamples to the Neggers-Stanley conjecture]
{Counterexamples to the Neggers-Stanley conjecture}
\author{Petter Br\"and\'en}
\address{Department of Mathematics, 
Chalmers University of Technology and G\"oteborg University,
S-412~96  G\"oteborg, 
Sweden}
\email{branden@math.chalmers.se}

\usepackage[english]{babel}
\usepackage{amsmath}
\usepackage{amsthm,amssymb,latexsym}
\usepackage{epic}
\usepackage{eepic}
\usepackage{xypic}
\usepackage{mathrsfs}

\newtheorem{conjecture}{Conjecture}
\newtheorem{theorem}{Theorem}
\theoremstyle{definition}
\newtheorem{example}{Example}
\newtheorem{remark}{Remark}

\newcommand{\CC}{\mathbb{C}}
\newcommand{\NN}{\mathbb{N}}
\newcommand{\e}{\ar@{-}}

\def\newop#1{\expandafter\def\csname #1\endcsname{\mathop{\rm
#1}\nolimits}}

\newop{des}
\newop{min}

\usepackage{graphicx}

\begin{document}

\maketitle
\nocite{*}\bibliographystyle{plain}

\begin{abstract}
The Neggers-Stanley conjecture (also known as the Poset conjecture) asserts 
that the polynomial counting the linear 
extensions of a partially ordered set on $\{1,2,\ldots,p\}$ by their number 
of descents has 
real zeros only. We provide counterexamples to this conjecture. 
\end{abstract}
\thispagestyle{empty}
\section{The Neggers-Stanley Conjecture}
Let $P$ be a poset (partially ordered set) on  
$[p]:=\{1,2,\ldots, p\}$. We will use the symbol $<$ to denote the 
usual order on the integers and the symbol $\prec$ to  
denote the partial order on $P$.  
The {\em Jordan-H\"older set}, $\mathcal{L}(P)$, of $P$ is the set of 
permutations   
$\pi$ of $[p]$   
which are linear extensions of $P$, i.e., if $i \prec j$ then 
$\pi^{-1}(i)<\pi^{-1}(j)$. A {\em descent} in a permutation 
$\pi$ is an index $i \in [p-1]$ such that 
$\pi(i)> \pi(i+1)$. Let $\des(\pi)$ denote the number of descents in 
$\pi$. The $W$-polynomial of $P$ is defined by 
$$
W(P,t) = \sum_{\pi \in \mathcal{L}(P)}t^{\des(\pi)}.
$$
The $W$-polynomials appear naturally in many combinatorial structures, 
see \cite{Brenti,Simion,Stanleythesis}, and are connected to Hilbert series 
of the Stanley-Reisner rings of simplicial complexes 
\cite[Section III.7]{Stanleygreen} and algebras with straightening laws 
\cite[Thm. 5.2.]{StanleyHilbert}.
\begin{example} Let $P_{2,2}$ be as in Fig. \ref{fig}. Then 
$$
\mathcal{L}(P_{2,2})=\{(1,3,2,4),(1,3,4,2),(3,1,2,4),(3,1,4,2),(3,4,1,2)\}
$$ 
so $W(P_{2,2},t)=4t+t^2$.
\end{example}
When $P$ has no relations, that is, when $P$ is the anti-chain on $[p]$, 
then $W(P,t)$ is the $p$th {\em Eulerian polynomial}. The Eulerian 
polynomials are known \cite{Harper} to have all zeros real. This is an 
instance where the Neggers-Stanley conjecture, also known as 
the Poset conjecture, holds: 
\begin{conjecture}[Neggers-Stanley]  
Let $P$ be a poset on $[p]$. Then all zeros of 
$W(P,t)$ are real.
\end{conjecture}
A poset $P$ is naturally labeled if 
$i\prec j$ implies $i<j$. 
The above conjecture was made by J. Neggers 
\cite{Neggers} in 1978 for naturally labeled posets and by R.P. Stanley in 
1986 for arbitrary posets on $[p]$. It has been proved for 
some special cases, see \cite{Brenti,Wagner}. 

If a polynomial $p(t)=a_0 +a_1t + \cdots + a_nt^n$ 
with non-negative 
coefficients is real-rooted it follows that the sequence 
$\{a_i\}_{i=0}^n$ is {\em unimodal} i.e., there is a $d$ such that 
$a_0 \leq a_1 \leq \cdots \leq a_{d-1} \leq a_d \geq a_{d+1} 
\geq \cdots \geq a_n$. This consequence of the Neggers-Stanley conjecture 
was recently proved \cite{ReinerWelker} (see also \cite{Branden}) for an 
important 
class of naturally labeled posets, namely {\em graded posets}. A 
naturally labeled poset $P$ is graded if all maximal chains in $P$ 
have the same length.      
\section{Counterexamples to the Neggers-Stanley conjecture}
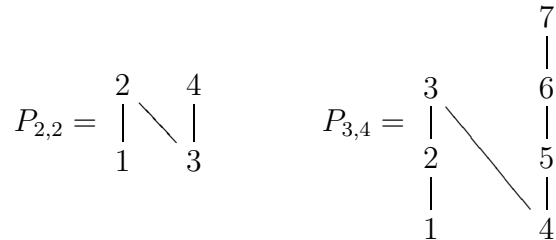
\begin{figure}\caption{\label{fig} The posets $P_{2,2}$ and $P_{3,4}$.}
\setlength{\unitlength}{10mm}
$$ 
P_{2,2}=\vcenter{\xymatrix@R=14pt@C=14pt{
              2& 4&&\\
              1\e[u]& 3\e[lu]\e[u]&&
          }}
P_{3,4}=\vcenter{\xymatrix@R=12pt@C=12pt{
                   && 7 \\              
              3    && 6\e[u]\\ 
              2\e[u]&&5\e[u]\\
              1\e[u]&& 4\e[lulu]\e[u]      
        }}
$$
\end{figure}

Let ${\bf m} \sqcup {\bf n}$ be the disjoint union of the chains  
$1\prec 2\prec \cdots \prec m$ and $m+1 \prec  n+2 \prec \cdots \prec m+n$. 
Simion \cite{Simion} proved that $W({\bf m} \sqcup {\bf n}, t)$ has 
real and simple zeros only. Let $P_{m,n}$ be the poset 
obtained by adding the relation $m+1 \prec m$ to the relations in 
${\bf m} \sqcup {\bf n}$, see 
Fig. \ref{fig}. 
The only linear extension 
of ${\bf m} \sqcup {\bf n}$ which is not a linear extension of 
$P_{m,n}$ is $(1,2, \ldots, m+n)$, which gives  
$$
W({\bf m} \sqcup {\bf n},t) =1 + W(P_{m,n},t).
$$
Let $\pi$ be a permutation. 
We say that $\pi(i)$ is a {\em descent top} and that $\pi(i+1)$ is a 
{\em descent bottom} if $i$ is a descent in $\pi$. A permutation in 
$\mathcal{L}({\bf m} \sqcup {\bf n})$ is uniquely determined by its descent 
tops (which 
are necessarily elements of $[m+n]\setminus[m]$) and its descent bottoms 
(which are elements of $[m]$). It follows that the number of 
permutations in $\mathcal{L}({\bf m} \sqcup {\bf n})$ with exactly $k$ 
descents is ${\binom m k}{\binom n k}$, so we have  
$$
W({\bf m} \sqcup {\bf n},t)= \sum_{k=0}^n {\binom m k}{\binom n k} t^k. 
$$ 
\begin{theorem}
Let $M>0$ be an integer. 
There is an integer $N=N(M)$ such that $W(P_{m,n},t)$ has more than 
$M$ non-real zeros 
whenever $\min(m,n)>N$.
\end{theorem}
\begin{proof}
Recall that the {\em Bessel function} of order $0$ is given by 
\begin{equation}\label{bessel}
J_0(z)= \frac 2 {\pi}\int_0^1\frac{\cos(zt)}{\sqrt{1-t^2}}dt
=\sum_{k=0}^{\infty}\frac 1{k!k!}\left(\frac{-z^2}{2}\right)^k.
\end{equation}
It follows from \eqref{bessel} that $J_0(z)$ has infinitely many  
zeros all of which are real and simple and that 
$|J_0(t)| \leq 1$ for all real $t$ with equality only if $t=0$. 
Hence 
the function 
$$
F(z)= \sum_{k=0}^{\infty}\frac{1}{k!k!}z^k.
$$
has infinitely many zeros all of which are negative and simple and 
$|F(t)| < 1$ for $t<0$.  

Let $f_{m,n}(t) = W(P_{m,n},t/mn)$. Then 
$$
f_{m,n}(t)= \sum_{k=1}^{\min(m,n)}\frac{\gamma_{m,k}\gamma_{n,k}}{k!k!}t^k,
$$
where $\gamma_{n,k}= (1-\frac 1 n)(1-\frac 2 n)\cdots (1-\frac{k-1}n)$. 
Let $\{(m_j,n_j)\}_{j=1}^{\infty} \subset \NN \times \NN$ be a sequence 
such that 
$\lim_{j \rightarrow \infty}\min(m_j,n_j)=\infty$. Then   
$$
\lim_{j \rightarrow \infty} f_{m_j,n_j}(z) = F(z)-1,
$$
where the sequence converges uniformly on all compact subsets of $\CC$.   
Let $(-a,0)$ be an interval containing more than $M$ zeros of $F(z)$. 
It follows from Hurwitz's theorem \cite[Thm. 1.3.8]{Rahman} 
that $f_{m_j,n_j}(z)+1$ has more than 
$M$ zeros in $(-a,0)$ for sufficiently large $j$ and by continuity that 
$|f_{m_j,n_j}(z)+1|<1$ for $z \in (-a,0)$. Thus by  
subtracting $1$ from $f_{m_j,n_j}(z)+1$ we will lose at least 
$M$ real zeros.  
\end{proof}
By applying Sturm's Theorem \cite[Section 10.5]{Rahman}
one can find specific counterexamples. The polynomial $W(P_{11,11},t)$ has two 
non-real zeros which are approximately 
$$
z= -0.10902 \pm 0.01308 i. 
$$
A counterexample with a polynomial of lower degree is  
$$
W(P_{36,6},t)=  216t + 9450t^2 + 142800t^3 + 
883575t^4 + 2261952t^5 + 1947792t^6.
$$
This polynomial has two non-real zeros. 
\begin{remark}
It should be noted that our counterexamples are not naturally labeled, since 
we have $n+1 \prec n$. The Neggers-Stanley conjecture is therefore still open 
for naturally 
labeled posets.
\end{remark}

\end{document}